\documentclass[a4paper,11pt]{article}
\usepackage[utf8]{inputenc}
\usepackage{amssymb,graphicx,epsfig,amsmath,amsfonts}
\usepackage{color}
\usepackage{colortbl}

\newcommand{\Rr}{{\mathbb{R}}}
\newcommand{\Nn}{{\mathbb{N}}}
\newcommand{\Cc}{{\mathbb{C}}}

\def\a2{{\alpha/2}}
\def\al{{\alpha}}
\def\bfc{{\mathbf{c}}}
\def\bfb{{\mathbf{b}}}
\def\bfq{{\mathbf{q}}}
\def\bfv{{\mathbf{v}}}
\def\bfz{{\mathbf{z}}}
\def\bfo{{\mathbf{0}}}
\def\bfE{{\mathbf{E}}}
\def\bfe{{\mathbf{e}}}
\def\a2{{\alpha/2}}
\def\hy{{\hat{y}}}
\def\GG{{\Gamma}}
\def\bfxi{\mbox{\boldmath$\xi$}}
\def\bftau{\mbox{\boldmath$\tau$}}
\def\qed{\hfill$\Box$}

\newtheorem{theorem}{Theorem}[section]
\newtheorem{remark}{Remark}[section]
\newtheorem{proposition}{Proposition}[section]

\title{A matrix method for fractional Sturm-Liouville problems on bounded domain.\thanks{{\em This work 
was supported by the GNCS-INdAM 2016 project ``Metodi numerici per operatori 
non-locali nella simulazione di fenomeni complessi''.}}}
\author{Paolo Ghelardoni\thanks{ Dipartimento di Matematica, Universit\`{a} di Pisa, Italy, paolo.ghelardoni@unipi.it}
\and 
Cecilia Magherini\thanks{Dipartimento di Matematica, Universit\`{a} di Pisa, Italy, cecilia.magherini@unipi.it.}}

\date{}

\begin{document}

\maketitle

\begin{abstract}
A matrix method for the solution of direct fractional Sturm-Liouville problems on bounded domain
is proposed where the fractional derivative is defined in the Riesz sense. The scheme is based 
on the application of the Galerkin spectral method of orthogonal polynomials. The order of convergence of the eigenvalue approximations with respect to the matrix size is studied. Some numerical
examples that confirm the theory and prove the competitiveness of the approach are finally 
presented.
\end{abstract}

\section{Introduction}
This paper concerns the numerical approximation of the eigenvalues of a time-independent one-dimensional
fractional Schroedinger equation defined on a bou\-nd\-ed interval which, without 
loss of generality, we assume to be $(-1,1).$  
This problem has several important applications. Among them we cite quantum mechanics with a Feynman path integral
over  L\'{e}vy trajectories, \cite{las1,las2}. Many other applications appear in mathematical
physics, biology and finance. \\

In more details, we shall consider the following eigenvalue problem
\begin{eqnarray}
 && \left(-\Delta\right)^{\a2} y(x) +q(x) y(x) = \lambda \,y(x), \qquad x\in D\equiv (-1,1), \label{fslp}\\
 &&\qquad \qquad  y(x) = 0 \quad \mbox{for each \,\,} x\in \Rr\setminus D, \label{bc}
\end{eqnarray}
where $\lambda$ and $y$ are an eigenvalue and a corresponding 
eigenfunction, respectively, $q$ represents the potential, and, for $\al \in (0,2),$ 
the {\em fractional Laplace operator} (or {\em quantum Riesz derivative})  is defined as
\begin{equation}\label{fdelta}
\left(-\Delta\right)^{\a2} y(x) \equiv \frac{1}{\eta(\al)} 
                   \lim_{\varepsilon\rightarrow 0^+} \int_{\Rr\setminus(-\varepsilon,\varepsilon)} \frac{y(x)-y(x-t)}{|t|^{1+\al}} dt,
\end{equation}
with
$$ \eta(\al) = -\frac{\pi^{1/2} \GG(-\a2)}{ 2^\al \GG((1+\al)/2)}
                = -2\GG(-\al) \cos(\al\pi/2).$$

Indeed several definitions of the fractional Laplacian can be found in the literature 
which are equivalent to (\ref{fdelta}) if $y$ and/or $\alpha$ verify suitable hypotheses 
(see, for instance, \cite{kwa15,luc13,MLP,Sam}).
One of them is given by the pseudo-differential operator with symbol $|\omega|^\alpha$, i.e.
$$
\left(-\Delta \right)^{\a2} y (x) = {\cal F}^{-1}\left\{ |\omega|^\alpha 
\hat{y}(\omega);x\right\} = \frac{1}{2\pi}\int_{-\infty}^{+\infty} |\omega|^\alpha e^{i\omega x}  
\hat{y}(\omega) \> d\omega, 
$$
where $\hat{y}$ is the Fourier transform of $y.$ It is known that this definition
is equivalent to (\ref{fdelta}) if $y \in L^s(\Rr)$ with $s\in[1,2],$ \cite{kwa15}.\\
Alternatively, if $\al \in (0,2)\setminus\left\{1\right\}$ then (\ref{fdelta}) can be written as
\begin{eqnarray} 
\left(-\Delta \right)^{\a2} y (x) = \frac{\sin(\al \pi /2)}{\sin(\al \pi)} \, \frac{d^n}{dx^n} I^{(n-\al)} y(x),  \label{laprie}
\end{eqnarray}
where $n=\lceil \al \rceil$ and
\begin{eqnarray*}
I^{(n-\al)} y(x) &=& \frac{1}{\Gamma(n-\al)} \int_{-\infty}^{+\infty} \left(\mbox{sign}(x-t)\right)^n |x-t|^{n-\al-1} y(t) dt\\
                  &=&  \int_{-\infty}^x  \frac{(x-t)^{n-\al-1}}{\Gamma(n-\al)} y(t) dt 
                    +  (-1)^n \int_x^{+\infty}  \frac{(t-x)^{n-\al-1}}{\Gamma(n-\al)} y(t) dt\\
                  &\equiv& I_{+}^{(n-\al)} y(x) + (-1)^n I_{-}^{(n-\al)}y(x).
\end{eqnarray*}
Here, $I_{\pm}^{(n-\al)}$ are the left- and 
right-sided {\em Riemann-Liouville fractional integrals} 
(sometimes called {\em Weyl integrals})
of order $n-\al.$\\
Justified by passage to the limit $\al \nearrow n,$ $I_{\pm}^{(0)}$ are defined as the identity
operator and, consequently, 
$\left(-\Delta\right)^{\a2}$ is set equal to $-\frac{d^2}{dx^2}$ for $\al=2,$ \cite{luc13}.
This implies that, for such value of $\al,$ (\ref{fslp})-(\ref{bc}) reduces
to the classical Sturm-Liouville problem in normal form with Dirichlet boundary conditions at both ends.\\

Concerning the special case of $q(x)\equiv 0$ in $D,$ sometimes referred to as the infinite 
potential well problem, it is known that the eigenvalues of (\ref{fslp})-(\ref{bc})
form an infinite sequence tending to infinity. More precisely, if we denote them with $\lambda_k$ then 
it is known that $0<\lambda_0<\lambda_1 \leq \lambda_2 \leq \lambda_3\leq \ldots$ and that the corresponding
eigenfunctions, say $y_k,$ form a complete orthonormal set in $L^2(D),$ \cite{BK04,BBKRSV}. 
Regarding the simplicity of the eigenvalues in \cite{kwa12}, see also the references therein, 
it was proved that this property is surely verified if $\al \in [1,2]$ and it was 
conjectured that indeed it holds for every $\al \in (0,2].$  Moreover, in the same paper, 
the following asymptotic law  
\begin{equation}\label{lamasq0}
\lambda_k = \left(\frac{(k+1)\pi}{2}- \frac{(2-\al)\pi}{8}\right)^\al + O\left(\frac{1}{k+1}\right)
\end{equation}
was determined (please observe that we number the eigenvalues starting from 
$k=0$ instead of $k=1$ as done in \cite{kwa12}).  It must be said that  an asymptotic growth like 
$((k+1)\pi/2)^\alpha$ was already proved in  \cite{CS05,DeB04}.\\
Now if $q\in L^2(D)$ then for the classical problem with $\al = 2$ it is known that  
\begin{equation}\label{lamasq}
\lambda_k(q) \approx \lambda_k(0) + \bar{q}, \qquad \bar{q} = \frac{1}{2} \int_{-1}^{1} q(x)dx, \qquad k\gg 0,
 \end{equation}
where $\lambda_k(q)$ and $\lambda_k(0)$ are the eigenvalues of index $k$ for the problems with potential $q$ and zero potential
in $D,$ respectively. More precisely, the residuals  $\delta_k=  \lambda_k(q) - \lambda_k(0) - \bar{q},$ $k \in \Nn_0,$ 
depend on $q-\bar{q}$ and constitute a square-summable sequence. In addition,
their rate of decrease is connected to the smoothness of $q$ over $[-1,1],$ \cite{PT}.
It is reasonable to assume that (\ref{lamasq}) holds true for each 
$\alpha \in(0,2]$ under the same hypothesis for $q.$\\

In the literature, the numerical schemes currently availables for the problem under consideration belong
to the family of so-called {\em matrix methods}, namely methods that discretize
the eigenvalue problem for the differential operator as an ordinary or a generalized 
matrix eigenvalue one. In particular, a number of finite difference schemes that
constitute a generalization of 
the classical three-point method (or discrete Laplacian in 1D) are available.
This is the case, for example, of the method proposed independently  by Ortigueira and by Zoia et.al. 
in \cite{Or06,ZRK07} and of the WSGD method (acronym for Weighted and Shifted Gr\"{u}nwald Difference) 
studied in \cite{TZD15}. In the former case, the discrete fractional Laplacian
is represented by the symmetric Toeplitz matrix $T_\al$ with symbol $(2(1-\cos\theta))^\a2$ 
(note that $2(1-\cos\theta)$
is the symbol associated to $T=\mbox{tridiag}(-1,2,-1)).$
The WSGD method, instead, provides an approximation of the left- or the right-sided Riemann-Liouville
fractional derivatives by using a suitable combination of the Gr\"{u}nwald and the 
shifted Gr\"{u}nwald difference schemes. 
In both the previous cases, in \cite{CD12,TZD15} it was proved that if 
$y$ is sufficiently regular over $\Rr$ 
and if $\alpha \in (1,2]$  then 
the error in the approximation of $\left(-\Delta\right)^\a2 y(x_n)$ behaves like $O(h^2),$
where $x_n$ and $h$ represent a meshpoint and the stepsize, respectively.
Unfortunately, the eigenfunctions of (\ref{fslp})-(\ref{bc}) are not smooth 
in proximity of the boundary of $D.$ For the problem with zero potential in $D,$ in fact,
it is known that there exist suitable constants $C_1,\,C_2,\,C_3$ and $\theta$ such that
\begin{equation} \label{yasint}
y(x) \approx \left\{\begin{array}{ll}
C_1 (1+x)^\a2 & \mbox{for } x \mbox{ close to } -1,\\
C_2 (1-x)^\a2 & \mbox{for } x \mbox{ close to } 1,\\
C_3 \sin(\lambda^{1/\al} x +\theta) & \mbox{for } x \in D \mbox{ away from the boundary,}
\end{array} \right. 
\end{equation}
see, for instance, \cite{kwa12} and \cite[Example~1]{kwa11}.
We expect that if $q\in L^2(D)$ and if $\|q-\bar{q}\|_2$ is not too large then the 
eigenfunctions have a similar behavior.
The lack of regularity near the boundaries of the domain is a peculiarity
of the solution of differential problems that involve fractional operators and,
in addition to the nonlocality of the latters, it represents a further important 
source of difficulties for their numerical treatment.
With reference to the  matrix methods previously mentioned, such 
behavior of $y$ determines an order reduction in the approximation of
its fractional Laplacian and consequently in the resulting numerical eigenvalues.
Alternative matrix methods are those proposed recently in \cite{BDM16,DZ15,GM14}.
In particular, the method in \cite{BDM16} is based on finite element 
approximations and it can be applied to problems in a generic dimension $d\ge 1,$
the approach considered in \cite{DZ15} is that of using suitable 
quadratures for the approximation of the integral in (\ref{fdelta})
and, finally, in \cite{GM14} a Control Volume Function approximation
with Radial Basis Function interpolation is proposed. All these schemes, however,
appear to be  of the first order, namely the error in the approximation of the eigenvalues
decreases like $N^{-1}$ where $N$ is the matrix size. \\

In this paper, we propose a matrix method based on  
the Galerkin spectral schemes named method of orthogonal polynomials 
in \cite{Pod} (see also the references therein).
Indeed, the principal idea has been recently presented in \cite{DKK} which
concerns the eigenvalue problem for the fractional Laplace operator in the unit ball 
of dimension $d\ge 1$ (so the potential is identically zero in $D$).
In such paper, it is proved that the eigenvalues provided by the matrix method 
with matrices of order $N,$ say $\lambda_k^{(N)},$ are such that $\lambda_k \leq \lambda_k^{(N)},$ for 
each $k< N.$ This is done by using the standard Rayleigh-Ritz variational method. 
In addition, the Aronszajn method of intermediate problems, see e.g. \cite{Bea}, is used for getting a lower bound 
for the eigenvalues. The aim pursued in \cite{DKK} is that of proving that if $1\leq d\leq 9$ and $\alpha=1$ or if
$1\leq d \leq 2$ and $\alpha\in(0,2]$ then the eigenfunctions corresponding to $\lambda_1$ are 
antisymmetric. In this paper, we are going to study such matrix method for $d=1$ from a numerical point of view.
More precisely, differently with respect to what has been done in \cite{DKK}, we shall consider a
generic potential $q\in L^2(D)$ and we will study the order of convergence of $|\lambda_k -\lambda_k^{(N)}|$
with respect to $N.$ Before proceeding, it must be said that the application of the method of
orthogonal polynomials to fractional eigenvalue problems has been recently considered also in \cite{ZK13,ZK15}
which however concern different fractional operators. 
For example, one of the generalization of 
the differential term $\frac{d}{dx}\left[p(x) \frac{d}{dx}\right]$ of the 
classical Sturm-Liouville operator considered in \cite{ZK13} (see also \cite{ka12}) is given by
$ \,_{~x}^{RL} D_1^{\a2} \left[ p(x) \,_{-1}^{~C} D_x^{\a2}\right].$
Here $_{~x}^{RL} D_1^{\a2}$ and $_{-1}^{~C} D_x^{\a2}$ are the right-sided  
Riemann-Liouville derivative of order $\a2,$ and the left-sided Caputo derivative of the same order, respectively.\\

The paper is organized as follows. In Section~\ref{polyorto} we introduce the approach based 
on the spectral  method of orthogonal polynomials. In Section~\ref{numerical} we derive the generalized 
matrix eigenvalue problem that discretize (\ref{fslp})-(\ref{bc}). Moreover, we describe how we have 
handled  a generic potential $q\in L^2(D)$ and we study the behavior of the entries 
in the resulting coefficient matrices. Section~\ref{error} is devoted to the analysis 
of the error in the eigenvalue approximations while Section~\ref{condizio} to the study 
of the conditioning of the numerical eigenvalues with respect to a perturbation of the potential.
Finally, in Section~\ref{example} we report the
results of several numerical examples that confirm the theory and prove the competitiveness of our
method.

\section{Spectral method of orthogonal polynomials}\label{polyorto}

By virtue of (\ref{yasint}), we consider the following expansion of an eigenfunction of (\ref{fslp})-(\ref{bc})
\begin{equation}\label{seriey}
 y(x) = (1-x^2)_+^{\a2} \sum_{n=0}^\infty c_{n} P_n^{(\a2,\a2)} (x).
\end{equation}
Here $a_+ = \max(a,0)$ and, for $\beta,\gamma>-1,$ $\left\{P_n^{(\beta,\gamma)}\right\}_{n\in\Nn_0}$ is the sequence 
of orthogonal Jacobi polynomials in $L^2(D,\omega)$  with weighting function 
$\omega(x) = (1-x)^\beta(1+x)^\gamma,$
i.e. $\langle P_n^{(\beta,\gamma)},P_m^{(\beta,\gamma)}\rangle_{\beta,\gamma} =0$ for each $n\neq m$ being
\begin{equation}\label{scalprod}
 \langle f,g\rangle_{\beta,\gamma} \equiv 
\int_{-1}^1 (1-x)^\beta (1+x)^\gamma f(x)\> g(x)\>dx, \quad  f,g\in L^2(D,\omega).
 \end{equation}
In particular, the following normalization 
$$ P_n^{(\beta,\gamma)}(1) = 1$$
will be used for such polynomials. If $\beta=\gamma$ then we shall use the simpler notation
\begin{equation}\label{notazione}
P_n^{(\beta)} \equiv P_n^{(\beta,\beta)}, \qquad \langle \cdot,\cdot\rangle_{\beta} \equiv\langle \cdot,\cdot\rangle_{\beta,\beta}.
\end{equation}
As we are going to show in Theorem~\ref{prop1}, the expansion in (\ref{seriey}) is favorable since
$$\left(-\Delta\right)^{\a2} \left((1-x^2)_+^{\a2} P_n^{(\a2)}(x) \right) 
 \propto P_n^{(\a2)}(x) \quad \mbox{for each  } x\in (-1,1).$$
Before this important result, for later convenience, we recall a list of known properties of the 
Jacobi polynomials revised according to the normalization that we have considered:  
\begin{enumerate}
 \item[\textbf{P1}:] $ \sigma_n \equiv \langle P_n^{(\a2)},P_n^{(\a2)}\rangle_{\a2}
 = \frac{2^{\al+1}\GG(n+1)\GG^2(\a2+1)}{(2n+\al+1) \GG(n+\al+1)};$
 \item [\textbf{P2}:] if $n$ is even then $P_n^{(\a2)}(x) = P_{n/2}^{(\a2,-1/2)}(2x^2-1);$
 \item [\textbf{P3}:] if $n$ is odd then $P_n^{(\a2)}(x) = xP_{(n-1)/2}^{(\a2,1/2)}(2x^2-1);$
 \item[\textbf{P4}:] for each $\beta,\gamma,\sigma>-1,$ \cite[16.4 formula (17)]{EMOT},
 \begin{eqnarray*}
 \nonumber \lefteqn{\hspace{-1cm}\langle P_r^{(\beta,\gamma)},P_s^{(\sigma,\gamma)}\rangle_{\beta+\sigma,\gamma} =}\\
   &&  \hspace{-1.25cm}\frac{(-1)^{r-s}\,2^{\gamma+\beta+\sigma+1} \,\GG(\beta+\sigma+1)\, \GG(\gamma+r+s+1) \GG(\beta+1) \GG(\sigma+1)}
  { \GG(\beta+r-s+1) \,\GG(\sigma+s-r+1) \,\GG(\gamma+\beta+\sigma+r+s+2)};
 \end{eqnarray*}
 \item[\textbf{P5}:] the polynomials $P_n^{(\a2)}$ verify the following recurrence relation
 \begin{eqnarray}
\nonumber P_{-1}^{(\a2)}(x) &\equiv& 0,\\
\label{P0}  P_0^{(\a2)}(x) &\equiv& 1,\\
\label{recurr} P_{n+1}^{(\a2)}(x) &=& \frac{2n+1+\al}{n+1+\al} x\, P_n^{(\a2)}(x) - \frac{n}{n+1+\al} P_{n-1}^{(\a2)}(x)\\
\label{recurr1}                    &\equiv& \zeta_{n,1} x P_n^{(\a2)}(x) - \zeta_{n,0}  P_{n-1}^{(\a2)}(x), \quad n\ge 0;
 \end{eqnarray}
 \item[\textbf{P6}:] $P_n^{(\a2)}(x) = \frac{\Gamma(\al+1)\Gamma(n+1)}{\Gamma(n+\al+1)}\,C_n^{(\a2+1/2)}(x),$
 where $C_n^{(\a2+1/2)}$ is the Gegenbauer polynomial of degree $n$ with its usual normalization, \cite{Andrews,Sz75};
 \item[\textbf{P7}:] $P_n^{(\beta,\gamma)}(x)$ coincides with the following 
  Gauss hypergeometric function
 \begin{equation}\label{PnHyp}
 P_n^{(\beta,\gamma)}(x) = \,_2F_1\left(-n,n+\beta+\gamma+1;\beta+1; (1-x)/2\right).
 \end{equation}
 \end{enumerate}
The latter property allows to extend the definition of $P_n^{(\beta,\gamma)}$ to all $\beta,\gamma \in \Rr$
with $-\beta \notin \Nn,$ \cite{Sz75}. \\

We can now prove the following \textit{spectral relationship} which is fundamental for 
the development of the method.

\begin{theorem}\label{prop1} If $\al\in (0,2]$  then for each $n\ge 0$ and each $x\in (-1,1)$
\begin{equation}\label{lapomP}
\left(-\Delta\right)^{\a2} \left((1-x^2)_+^{\a2} P_n^{(\a2)}(x) \right) 
 = \mu_n P_n^{(\a2)}(x),
\end{equation}
where
\begin{equation}\label{mun}
\mu_n = \frac{\GG(n+\al+1)}{\GG(n+1)}.
\end{equation}
\end{theorem}
\underline{Proof} 
 If $\al=2$ then $P_n^{(1)} (x) = L_{n+1}'(x)$ where $L_{n+1}(x)$ is the Legendre polynomial of 
 degree $n+1,$ with a suitable normalization, and $\mu_n = (n+2)(n+1).$ 
 It follows that (\ref{lapomP})-(\ref{mun}) reduce to the well-known identity
 $$-\frac{d^2}{dx^2} \left((1-x^2) L_{n+1}'(x)  \right) = (n+2)(n+1) L_{n+1}'(x), \qquad x\in (-1,1).$$
Let's consider the case $\al\in (0,2)$ and $\al\neq 1.$ After some computations, by using 
(\ref{laprie}) and (\ref{PnHyp}), one obtains that (\ref{lapomP})-(\ref{mun}) are an 
application of Theorems~6.2 and 6.3 in \cite{Pod}. The special value $\al=1$ follows by continuity.\\
Alternatively, by virtue of properties \textbf{P2} and \textbf{P3},  the statement 
is an application of Theorem~3 in \cite{DKK_1} for every $\al \in (0,2].$\qed\\

Now, if $y$ satisfies (\ref{fslp})-(\ref{bc}) with eigenvalue $\lambda$ then for each $m \in \Nn_0$ 
\begin{eqnarray}\label{debole1}
 \langle P_m^{(\a2)},\left(-\Delta\right)^{\a2} y\rangle_{\a2} + \langle P_m^{(\a2)},q y\rangle_{\a2} 
  = \lambda \langle P_m^{(\a2)},y\rangle_{\a2},
\end{eqnarray}
see (\ref{scalprod})-(\ref{notazione}). 
Therefore, from (\ref{seriey}) and (\ref{lapomP})-(\ref{mun})
one gets that the first term in the previous equation reduces to 
$$\langle P_m^{(\a2)},\left(-\Delta\right)^{\a2} y\rangle_{\a2} = \mu_m \langle P_m^{(\a2)},P_m^{(\a2)}\rangle_{\a2}c_m\equiv a_m c_m$$
where, see property {\bf P1}, 
\begin{equation}\label{amm1}
a_m = \mu_m \langle P_m^{(\a2)},P_m^{(\a2)}\rangle_{\a2} = \frac{2^{\al+1}\GG^2(\a2+1)}{2m+\al+1}.
\end{equation}

Concerning the inner product on the right-hand side of (\ref{debole1}), from (\ref{seriey}) it follows that
\begin{equation}\label{bmn}
\langle P_m^{(\a2)},y\rangle_{\a2} \equiv \sum_{n=0}^\infty b_{mn} c_{n},\qquad   b_{mn}\equiv \langle P_m^{(\a2)},P_n^{(\a2)}\rangle_{\al} .
\end{equation}
It is evident that $b_{mn}=0$ for each $m$ and $n$ such that $m+n$ is odd. 
Moreover, the application of properties {\bf P2}--{\bf P4} allows to determine the remaining values analytically.\\
Let us consider the case where
$m=2r$ and $n=2s.$ By using property {\bf P2}, one deduces that
\begin{eqnarray}
\nonumber \hspace{-.7cm}\lefteqn{b_{mn} = \int_{-1}^{1} (1-x^2)^{\al}P_m^{(\a2)}(x)P_n^{(\a2)}(x) dx}\\
\nonumber &&\hspace{-.9cm} = 2\int_{0}^{1} (1-x^2)^{\al}P_r^{(\a2,-1/2)}(2x^2-1)P_s^{(\a2,-1/2)}(2x^2-1) dx\\
\nonumber     &&\hspace{-.9cm} = \frac{1}{2^{1/2+\al}} \int_{-1}^1 (1-t)^{\al} (1+t)^{-1/2} P_r^{(\a2,-1/2)}(t)P_s^{(\a2,-1/2)}(t) dt\\
\label{bmn_pari} &&\hspace{-.9cm} = \frac{(-1)^{(m-n)/2}\, \GG(\al+1)\, \GG((m+n+1)/2)\,\GG^2(\a2+1)}
              {\GG((\al-m+n)/2+1)\,\GG((\al-n+m)/2+1)\GG((m+n+3)/2+\al)} .
\end{eqnarray}  

\noindent In particular, the last equality follows from property {\bf P4} with $\gamma=-1/2,$ $\beta=\sigma=\a2,$
$r=m/2,$ and $s=n/2.$ It must be said that the previous formula was already determined in \cite{DKK_1}, with suitable changes 
in the notation and by considering the different normalization of the Jacobi polynomials.\\ 
Now (\ref{bmn_pari}) holds true also in the case where $m$ and $n$ are odd, i.e. $m=2r+1$ and $n=2s+1.$ 
In fact, from property {\bf P3}, we get
\begin{eqnarray*}
b_{mn} &=& 2\int_{0}^{1} (1-x^2)^{\al} x^2 P_r^{(\a2,1/2)}(2x^2-1)P_s^{(\a2,1/2)}(2x^2-1) dx\\
 &=& \frac{1}{2^{3/2+\al}} \int_{-1}^1 (1-t)^{\al} (1+t)^{1/2} P_r^{(\a2,1/2)}(t)P_s^{(\a2,1/2)}(t) dt
\end{eqnarray*}
which one can verify to be equal to the right-hand side of (\ref{bmn_pari}) by using property {\bf P4} 
with $\gamma=1/2,$ $\beta=\sigma=\a2,$ $r=(m-1)/2,$ and $s=(n-1)/2.$ \\

We observe that the application of Euler's reflection formula $\GG(1-z)\GG(z) = \pi/\sin(\pi z),$ 
see also \cite[eq.(4.21)-(4-23)]{Or06}, allows to get that  if $m+n$ is even then
$$ \frac{(-1)^{(m-n)/2}}{\GG((\al-m+n)/2+1)\,\GG((\al-n+m)/2+1)} = 
-\frac{\sin(\pi\a2)}{\pi} \, \frac{\GG\left(\left|\frac{n-m}2\right| -\frac{\al}2\right)}{\GG\left(\left|\frac{n-m}2\right| +\frac{\al}2+1\right)}.
$$
This implies that the coefficient $b_{mn}$ in (\ref{bmn_pari}) can be written as 
\begin{equation}\label{bmn_alt}
b_{mn} = \theta_\al h_{m+n} t_{|n-m|} 
 \end{equation}
where
\begin{eqnarray}
\nonumber
\theta_\al &=& -\frac{\sin(\pi\a2) \GG(\al+1)\,\GG^2(\a2+1)}{\pi}, \\
\label{hmn}h_{m+n} &=& \frac{ \GG((m+n+1)/2)}{\GG((m+n+3)/2+\al)}, \\
\label{tmn} t_{|n-m|} &=&  \frac{\GG\left(\left|\frac{n-m}2\right| -\frac{\al}2\right)}
              {\GG\left(\left|\frac{n-m}2\right| +\frac{\al}2+1\right) }.
\end{eqnarray}

\section{Numerical scheme}\label{numerical}
In order to get a numerical method for the approximation of the eigenvalues and of the eigenfunctions of the 
fractional Sturm-Liouville problem, we truncate the series in (\ref{seriey}), i.e. we 
look for an approximation of $y$ of the form
\begin{equation}\label{yN}
y(x) \approx y^{(N)}(x) = (1-x^2)_+^{\a2}\sum_{n=0}^{N-1} \xi_{n,N} P_n^{(\a2)}(x)
\end{equation}
where the coefficients $\xi_{n,N}$ are determined by imposing that (\ref{debole1}) holds true for
$m=0,\ldots,N-1$ with $y$ and $\lambda$ replaced by $y^{(N)}$ and $\lambda^{(N)}$ respectively. 
This leads to a generalized matrix eigenvalue problem of the form
\begin{equation}\label{genmat}
\left(A_N + Q_N\right) \bfxi_N = \lambda^{(N)} B_N \bfxi_N,
\end{equation}
where $\bfxi_N = \left(\xi_{0,N},\xi_{1,N},\ldots,\xi_{N-1,N}\right)^T$ and, see (\ref{amm1})-(\ref{bmn}), 
$$A_N = \mbox{diag}\left(a_0,\ldots,a_{N-1}\right), \qquad 
B_N = \left(b_{mn}\right)_{m,n=0,\ldots,N-1}.$$
Finally, the entries of $Q_N$ are given by
\begin{equation} \label{qmn} 
q_{mn} = \langle P_m^{(\a2)}, q P_n^{(\a2)}\rangle_{\al} = \int_{-1}^1 (1-x^2)^{\al} q(x) P_m^{(\a2)}(x) P_n^{(\a2)}(x) dx.
\end{equation}
Clearly, they are not known in closed form for a generic potential $q(x).$ 
We will talk about their approximation in Subsection~\ref{ssQ}.

\begin{remark}
$B_N$ is permutation similar to a $2 \times 2$ block diagonal matrix. The same holds true 
for $Q_N$ if the potential is an even function.
\end{remark}

\begin{remark}\label{Bsdp}
$B_N$ is symmetric positive definite since  
$$ \bfv^{T}B_N\bfv = \int_{-1}^{1} (1-x^2)^{\al} v^2(x) dx >0, \qquad v(x) = \sum_{n=0}^{N-1} v_n P_n^{(\a2)}(x),$$
for each $\bfv=\left(v_0,\ldots,v_{N-1}\right)^T \in \Rr^{N}\setminus \left\{\bfo_N\right\}$
and its simmetry is obvious.
\end{remark}

We observe that from (\ref{bmn_alt})--(\ref{tmn}) it is not difficult to deduce that $B_N$ is an Hadamard product between 
an Hankel matrix and a symmetric Toeplitz one. Moreover, its nonzero entries can be computed 
with a computational cost rather low by using the following recurrence relations 
\begin{eqnarray*}
h_{m+n} &=& \frac{m+n-1}{m+n+1+2\al} \,h_{m+n-2},\\
t_{|n-m|} &=& \frac{|n-m|-\al-2}{|n-m|+\al} \, t_{|n-m|-2},
\end{eqnarray*}
that, in addition, allow to avoid problems of overflow and/or underflow.\\
Finally, for the error analysis in the eigenvalue approximations, 
it is important to analyze the behavior of such coefficients and, consequently, of $b_{mn}$ 
when $m+n$ and/or $|m-n|$ become large. 
We recall the following expansion of the ratio of two gamma functions  
$$\frac{\GG(z+a)}{\GG(z+b)} = z^{a-b} \left(1 + \frac{(a-b)(a+b-1)}{2z} + O (|z|^{-2})\right),\quad
z\neq 0.$$
Its application to (\ref{bmn_alt}) and (\ref{hmn})-(\ref{tmn}), for $m+n$ even, allows to obtain that 
\begin{itemize}
 \item if $m+n>0$ then 
 $$h_{m+n} = 2^{\al+1}(m+n)^{-\al-1} \left(1+O((m+n)^{-1})\right);$$
\item if $|n-m|>0$ then 
$$ t_{|n-m|} = 2^{\al+1}|n-m|^{-\al-1} \left(1 + O((n-m)^{-2})\right);$$
\item if $n> m$ then
\begin{equation}\label{bmn_as}
b_{mn} = b_{nm} = \theta_{\al} \left(\frac{4}{n^2-m^2}\right)^{\al+1}\left(1+O((m+n)^{-1}) + O((n-m)^{-2})\right).
\end{equation}
\end{itemize}

\subsection{Computation and properties of the entries of $Q_N$}\label{ssQ}
We are now going to talk about possible techniques for computing the entries of the matrix $Q_N$ 
in (\ref{genmat}) and about their asymptotic behavior.\\
Considering the definition of $q_{mn}$ in (\ref{qmn}), the first idea for its approximation is 
trivially that of applying a Jacobi quadrature rule with
weighting function $(1-x^2)^{\al}.$ This is surely a possibility which, however, 
requires the application of a formula with degree of precision  
rather large since the integrand is $q(x) P_m^{(\a2)}(x)P_n^{(\a2)}(x).$ \\
A second approach is suggested by the following results, \cite{Golub}. 
\begin{proposition} For each $m,n\ge 0,$ let $q_{mn}$ be defined as in (\ref{qmn}) and 
$q_{m,-1}=q_{-1,m}=0.$ Then, see (\ref{recurr})-(\ref{recurr1}), 
\begin{equation} \label{recqmn}
q_{m,n+1} = \frac{\zeta_{n,1}}{\zeta_{m,1}} q_{m+1,n} + \frac{\zeta_{n,1}\zeta_{m,0}}{\zeta_{m,1}} q_{m-1,n}
-\zeta_{n,0} q_{m,n-1}, \quad m,n\ge 0.
\end{equation}
\end{proposition}
\underline{Proof} 
From (\ref{qmn}) and (\ref{recurr})-(\ref{recurr1}), we get 
\begin{eqnarray*}
q_{m,n+1} &=& \langle q \,P_m^{(\a2)}, P_{n+1}^{(\a2)} \rangle_{\al} \\
&=& \zeta_{n,1}  \langle q \,P_m^{(\a2)}, x\,P_{n}^{(\a2)} \rangle_{\al}
-\zeta_{n,0}  \langle q \,P_m^{(\a2)}, P_{n-1}^{(\a2)} \rangle_{\al}\\
&=& \zeta_{n,1}  \langle q \,x\,P_m^{(\a2)}, P_{n}^{(\a2)} \rangle_{\al}
-\zeta_{n,0}  q_{m,n-1}.
\end{eqnarray*}
Therefore, the statement is a consequence of the fact that 
$$ x P_m^{(\a2)}(x) = \frac{1}{\zeta_{m,1}} \left( P_{m+1}^{(\a2)}(x) + \zeta_{m,0} P_{m-1}^{(\a2)}(x) \right). $$
\qed\\

\begin{proposition}\label{propQinf} Let $Q_\infty = \left(q_{mn}\right)_{m,n \in \Nn_0},$ 
$B_\infty = \left(b_{mn}\right)_{m,n \in \Nn_0}$ and  
\begin{equation} \label{bfqb}
\bfq_n \equiv \left(\begin{array}{c}q_{0n}\\q_{1n}\\ \vdots\end{array}\right)\in \ell_\infty, \qquad 
\bfb_n \equiv \left(\begin{array}{c}b_{0n}\\b_{1n}\\ \vdots\end{array}\right)\in \ell_\infty,
\end{equation}
i.e. let $\bfq_n$ and $\bfb_n$ be the $n$-th column of $Q_\infty$ and $B_\infty,$ respectively.
If we define the following linear tridiagonal operator $ \bfz \in \ell_\infty \mapsto \mathcal{H}\bfz \in \ell_\infty$ where,
see (\ref{recurr})-(\ref{recurr1}),
\begin{eqnarray} \nonumber
\mathcal{H} &=& \left( \begin{array}{cccccc}
            0 & \frac{1} {\zeta_{0,1}}\\
            \frac{\zeta_{1,0}} {\zeta_{1,1}} & 0 & \frac{1} {\zeta_{1,1}} \\
            &\frac{\zeta_{2,0}} {\zeta_{2,1}} & 0 & \frac{1} {\zeta_{2,1}} \\
            && \ddots & \ddots & \ddots\\
            &&& \ddots & \ddots & \ddots\\
            \end{array}\right)\\      
&=&
\left( \begin{array}{cccccc}
            0 & 1\\
            \frac{1} {3+\al} & 0 & \frac{2+\al} {3+\al}\\
            &\frac{2} {5+\al} & 0 & \frac{3+\al} {5+\al}\\
            && \ddots & \ddots & \ddots\\
            &&& \ddots & \ddots & \ddots\\
            \end{array}\right)   \label{H},
\end{eqnarray}
then
\begin{enumerate}
 \item $\mathcal{H} \,\bfE = \bfE$ where $\bfE=(1,1,\ldots)^T;$
 \item by setting $\bfq_{-1} = (0,0,\ldots)^T,$ one gets 
 \begin{equation}\label{recqn}
\bfq_{n+1} = \zeta_{n,1} \mathcal{H}\bfq_n -\zeta_{n,0} \bfq_{n-1},\qquad \mbox{for each  } n\ge 0,  
 \end{equation}
 \item $\bfq_n =  P_n^{(\a2)} (\mathcal{H}) \,\bfq_0,$
 \item $\bfb_n =  P_n^{(\a2)} (\mathcal{H}) \,\bfb_0.$

\end{enumerate}
\end{proposition}
\underline{Proof} 
 The first result is trivial while the second and, consequently, the third ones follow from 
(\ref{recqmn}).
Concerning the last point, it is sufficient to observe that $B_\infty=Q_\infty$ if $q(x) \equiv 1.$ \qed\\

Let now assume that we know the first $S+1$ entries in $\bfq_0.$ By using (\ref{recqn}) with $n=0$ we can compute
the first $S$ entries of $\bfq_1.$ At this point, from the same formula with $n=1,$
we determine the values of the first $S-1$ entries of $\bfq_2$ and so on. By  observing that $Q_N$ in (\ref{genmat}) 
is of order $N,$ we deduce that it is entirely determined once the values of $q_{m0}$ for $m\leq S=2N$ are known. 
Clearly, in the actual implementation, the symmetry of $Q_N$ is taken into account.\\

The following result will be useful for the error analysis in the eigenvalue approximation.

\begin{proposition}\label{Propqas} Suppose that $q$ is analytic inside and on the Bernstein ellipse 
${\cal{E}}_\rho$ given by
\[ {\cal{E}}_\rho =\left\{ z\in \Cc\> |\> z=\frac{1}{2}(\rho e^{i\theta}+\rho^{-1}e^{-i\theta}),
\> 0\leq \theta\leq 2\pi \right\} \> , \]
with $\rho>1.$ If $m+n\gg 1$ and $m\neq n$ then $q_{mn}$ defined in (\ref{qmn}) satisfies
 \begin{equation}\label{qmn_as}
 q_{mn} = O\left( (m+n)^{-\al-1} \left(|m-n|\right)^{-\al-1}\right). 
 \end{equation}
 \end{proposition}
\underline{Proof} 
The regularity of $q$ implies that its Fourier-Jacobi expansion
\[ q(x)= \sum_{\ell =0}^\infty \gamma_\ell P_\ell^{(\alpha/2)}(x)\> , \quad \gamma_\ell
=\frac{<q,P_\ell^{(\alpha/2)}>_\a2}{<P_\ell^{(\a2)},P_\ell^{(\a2)}>_\a2},\]
converges in uniform norm with an exponential decay of the coefficients $\gamma_\ell.$
More precisely, see \textbf{P6}, in \cite{Wang}
it is proved that
\begin{equation}\label{gamas}
 \gamma_\ell \sim \ell^{(\alpha +3)/2} \rho^{-(\ell +1)} \> .
\end{equation}
Now, from the definition of the entries in $\bfq_0$ and $\bfb_\ell$ and from the fourth point in Proposition~\ref{propQinf} we get 
\begin{equation}\label{qH}
\bfq_0 =  \sum_{\ell =0}^\infty \gamma_\ell \bfb_\ell = 
             \left(\sum_{\ell =0}^\infty \gamma_\ell P_\ell^{(\a2)} (\mathcal{H})\right) \bfb_0 \equiv q(\mathcal{H}) \bfb_0.
\end{equation}     
In addition, by applying the principle of induction and by using (\ref{P0})--(\ref{recurr1}) and (\ref{H}) in a way similar to what was done in the proofs 
of the previous two propositions, one obtains that the entries of $P_\ell^{(\a2)} (\mathcal{H})$ are given by
\[ \left(P_\ell^{(\a2)} (\mathcal{H})\right)_{mj} = 
\frac{<P_\ell^{(\a2)}P_m^{(\a2)},P_j^{(\a2)}>_\a2}{<P_j^{(\a2)},P_j^{(\a2)}>_\a2}, \qquad m,j \in \Nn_0. \]
Clearly they are nonzero if and only if $m$ and $j$ are such that $m+\ell+j$ is even and the sum of any two of $\ell,m,j$ is 
not less than the third. Moreover, these entries are known in closed form thanks to properties {\bf P1}, {\bf P6} and to \cite[Corollary 6.84, pag. 321]{Andrews}. 
In particular, long and tedious computations allow to get that if $k\in\left\{-\ell,-\ell+2,\ldots,\ell\right\}$ then
$$
\left(P_\ell^{(\a2)} (\mathcal{H})\right)_{j-k,j} = \Upsilon_{\ell,k} \left(1 + \frac{(\al+1) k}{2j} + O(j^{-2})\right),  
\quad  j \gg \ell,
$$
where
$$
\Upsilon_{\ell,k}  = \frac{{\cal B}\left(\frac{\ell+k+\al+1}{2},\frac{\ell-k+\al+1}{2}\right)}
{(\ell+1) {\cal B}\left(\frac{\al+1}{2},\frac{\al+1}{2}\right) {\cal B}\left(\frac{\ell+k+2}{2},\frac{\ell-k+2}{2}\right)}.
$$
being ${\cal B}(\cdot,\cdot)$ the beta function. This implies that, see (\ref{qH}),  
$$ \lim_{j\rightarrow \infty}  \left(q (\mathcal{H})\right)_{j-k,j} = \sum_{\ell=k,k+2,\ldots} \gamma_\ell \Upsilon_{\ell,k}.$$
Specifically,  thanks to (\ref{gamas}), the modulus of the entries of $q(\cal{H})$ decays exponentially when going away from the main 
diagonal, i.e. when $|k|$ increases.
Therefore, by considering that from the third and the fourth points in Proposition~\ref{propQinf} one gets
$$
\bfq_n =  P_n^{(\a2)} (\mathcal{H}) \,\bfq_0 =  P_n^{(\a2)} (\mathcal{H}) q(\mathcal{H}) \bfb_0
= q(\mathcal{H})\left( P_n^{(\a2)} (\mathcal{H}) \bfb_0\right) = q(\mathcal{H}) \bfb_n,
$$
the estimate in (\ref{qmn_as}) follows from (\ref{bmn_as}).\qed\\

\begin{remark} \label{qpoli}
It is important to underline the fact that if $q$ is a polynomial then (\ref{recqn}) and (\ref{qH}) allow to 
compute $Q_N$ in a simple way.
\end{remark}

\section{Error analysis}\label{error}
This section is devoted to the analysis of the order of convergence of $\lambda^{(N)}$ versus $\lambda$
as $N$ increases. We will always suppose that $q$ satisfies the following hypotheses:
\begin{itemize} 
 \item[{\bf H1}] $q$ is analytic inside and on the Bernstein ellipse ${\cal E}_\rho$ with $\rho >1$ (see Proposition~\ref{Propqas});
 \item [{\bf H2}] $\|q-\bar{q}\|_2$ is not too large where $\bar{q}$ is the mean value of $q$ defined in (\ref{lamasq}).
\end{itemize}

Let $\bfc_N=\left(c_0,c_1,\ldots,c_{N-1}\right)^T$ be the vector containing the first $N$ coefficients of
the expansion in (\ref{seriey}) of the exact eigenfunction corresponding to $\lambda$ and, see (\ref{genmat}), let
\begin{equation}\label{tau_def}
\bftau^{(N)} = \left(A_N+Q_N\right) \bfc_N -\lambda B_N \bfc_N
\end{equation}
be the local truncation error.
By applying standard arguments one obtains
\begin{equation}\label{errlam}
|\lambda-\lambda^{(N)}| = \left|\left(\bfxi_N^T \bftau^{(N)}\right)/\left(\bfxi_N^T \,B_N\,\bfc_N\right)\right|\> .
\end{equation}
Clearly, the error in the eigenvalue approximation is independent of the
normalization considered for its exact and numerical eigenfunctions. Therefore, for later convenience, 
we normalize $y(x)$ as follows 
\begin{equation}\label{normalcn}
y(x) = (1-x^2)^\a2 \,\hat{y}(x), \quad \hat{y}(1) = 1 \qquad  \Longleftrightarrow \qquad 
\sum_{n=0}^\infty c_n =1 
\end{equation}
and, ideally, we scale the numerical eigenfunction so that 
$$\xi_{s,N} = c_s \qquad \mbox{ where $s$ is such that} \quad c_0=\ldots=c_{s-1}=0, \quad c_s\neq 0.$$
Let us consider, first of all, the denominator in (\ref{errlam}). From Remark~\ref{Bsdp},
it follows that 
\begin{eqnarray}
\nonumber \lim_{N\rightarrow\infty} \bfxi_N^T B_N \bfc_N &=& 
 \lim_{N\rightarrow\infty} \left(\bfc_N^T B_N \bfc_N + (\bfxi_N-\bfc_N)^T B_N \bfc_N \right)\\
\label{limden}  &=& \|y\|_2^2 + \lim_{N\rightarrow\infty} (\bfxi_N-\bfc_N)^T B_N \bfc_N = \|y\|_2^2 
 \end{eqnarray}
provided that, by Cauchy-Schwarz, the $L_2$-norm of  
$(1-x^2)^{\alpha/2}\sum_{n=0}^{N-1} (\xi_{n,N}-c_n)P_n^{(\alpha/2)}(x)$ approaches zero 
(at least slowly) as $N$ tends to infinity. This implies that the denominator in (\ref{errlam}) is kept 
away from zero as $N$ increases. \\
Concerning the local truncation error, from (\ref{seriey}), (\ref{debole1})--(\ref{bmn}) and (\ref{qmn})
one gets 
$$ a_mc_m+\sum_{n=0}^\infty q_{mn} c_n=\lambda\sum_{n=0}^\infty  b_{mn} c_n\> .$$
It follows that the entries of $\bftau^{(N)}$ in (\ref{tau_def}) can be written as
\begin{equation}\label{tau_m2}
\tau_m^{(N)} \equiv \bfe_{m+1}^T \bftau^{(N)} =\sum_{n=N}^\infty (\lambda b_{mn} - q_{mn} ) c_n \equiv \sum_{n=N}^\infty u_{mn} c_n , \>\>\> m=0,..., N-1.
\end{equation}
 
We have already studied the behaviors of $b_{mn}$ and of $q_{mn},$ 
see (\ref{bmn_as})--(\ref{qmn_as}). 
It remains to establish how  $|c_n|$ decreases as $n$ increases.
Let us separate the even and the odd parts of $\hat{y}(x)$, $y(x)$  and of the potential
\begin{eqnarray*}
\hat{y}(x) &=& \hat{y}_e(x) +\hat{y}_o(x), \qquad \hat{y}_e(x) = (\hat{y}(x)+\hat{y}(-x))/2, \\ 
y(x) &=& y_e(x) +y_o(x), \qquad y_e(x) = (1-x^2)^\a2 \hat{y}_e(x), \\
q(x) &=& q_e(x) + q_o(x), \qquad q_e(x) = (q(x)+q(-x))/2. 
\end{eqnarray*}
It is not difficult to verify that
\begin{eqnarray}\label{probpari}
(-\Delta)^\a2  y_e(x) &=& (\lambda-q_e(x)) y_e(x) - q_o(x) y_o(x), \\
(-\Delta)^\a2 y_o(x) &=& (\lambda-q_e(x)) y_o(x) - q_o(x) y_e(x). \label{probdisp}
\end{eqnarray}
Now, if $q$ satisfies {\bf H1} and {\bf H2} then there exist suitable 
$\beta_e,\beta_o>0$ such that when $x \rightarrow \pm 1^{\mp}$ it results 
\begin{eqnarray*}
\hy_e(x) &=& \hy_e(1)P_0^{(\a2)}(x) + O((1-x^2)^{\beta_e}), \\
\hy_o(x) &=& \hy_o(1)P_1^{(\a2)}(x) + O(x(1-x^2)^{\beta_o}),\\
q_e(x) &=& q_e(1)P_0^{(\a2)}(x) + O((1-x^2)), \\
q_o(x) &=& q_o(1)P_1^{(\a2)}(x) + O(x(1-x^2)). 
\end{eqnarray*}
In addition, from the recurrence relation in (\ref{recurr}), one gets
$$ \left(P_1^{(\a2)}(x)\right)^2 = x P_1^{(\a2)}(x) = \frac{2+\al}{3+\al} P_2^{(\a2)} (x) + \frac{1}{3+\al} P_0^{(\a2)}(x).$$
It follows that in proximity of the extremes  
the terms on the right-hand side of (\ref{probpari})-(\ref{probdisp}) can be written respectively as
\begin{eqnarray}
\nonumber \lefteqn{(\lambda - q_e(x)) y_e(x) - q_o(x) y_o(x)}\\ 
&&\qquad \qquad =  (1-x^2)^\a2 \left(\nu_0 P_0^{(\a2)}(x) -\nu_2 P_2^{(\a2)}(x)\right) + r_e(x)\> , \label{rhsp}\\
\lefteqn{(\lambda - q_e(x)) y_o(x) - q_o(x) y_e(x) = (1-x^2)^\a2 \nu_1 P_1^{(\a2)}(x) + r_o(x)\> ,} \label{rhsd} 
\end{eqnarray}
where 
\begin{eqnarray}
\label{nu0a} \nu_0 &=& (\lambda - q_e(1))\hy_e(1) - q_o(1)\hy_o(1)/(3+\al),\\
\label{nu1b} \nu_1 &=& (\lambda - q_e(1))\hy_o(1) - q_o(1)\hy_e(1),\\
\nonumber \nu_2 &=& q_o(1)\hy_o(1) (2+\al)/(3+\al),
\end{eqnarray}
and $r_e(x)$ and $r_o(x)$ approach zero faster than $(1-x^2)^\a2$ as $x\rightarrow \pm 1^{\mp}.$
All these arguments allow to prove the following result.

\begin{theorem}\label{thcnas}
If $\al \in (0,2)$ and if $q$ is such that {\bf H1} and {\bf H2} hold true then
\begin{equation} \label{cn_as}
c_n = O(n^{-2\al -1}), \quad  n\gg k,
\end{equation}
where $k$ is the index of the eigenvalue.
\end{theorem}
\underline{Proof} 
First of all, we observe that for each $i$ and each $x\in (-1,1)$ 
\begin{eqnarray*}
(1-x^2)^\a2 P_i^{(\a2)} (x) &=& \sum_{n=0}^\infty \frac{\langle P_n^{(\a2)}, (1-x^2)^\a2 P_i^{(\a2)} \rangle_\a2}{\langle P_n^{(\a2)},P_n^{(\a2)}\rangle_\a2} \,P_n^{(\a2)}(x)\\
                            &=& \sum_{n=0}^\infty \frac{\langle P_n^{(\a2)}, P_i^{(\a2)} \rangle_\al}{\langle P_n^{(\a2)},P_n^{(\a2)}\rangle_\a2} \,P_n^{(\a2)}(x)\\
                            &=& \sum_{n=0}^\infty \frac{b_{ni}}{\sigma_n} P_{n}^{(\a2)}(x),
\end{eqnarray*}
see (\ref{bmn}) and property {\bf P1}. As a consequence, by inserting (\ref{rhsp}) and 
(\ref{rhsd}) into (\ref{probpari}) and (\ref{probdisp}), respectively,
and by considering the expansion of both sides of the equations in terms of the Jacobi 
polynomials, we obtain that for $x\rightarrow \pm 1^\mp$ 
\begin{eqnarray*}
\sum_{n=0}^\infty \mu_{2n}c_{2n} P_{2n}^{(\a2)} (x) &=& \sum_{n=0}^\infty  \left(\frac{\nu_0 b_{2n,0}- \nu_2 b_{2n,2}}{\sigma_{2n}} + r_{2n}\right) P_{2n}^{(\a2)}(x),\\
\sum_{n=0}^\infty \mu_{2n+1}c_{2n+1} P_{2n+1}^{(\a2)} (x) &=& \sum_{n=0}^\infty  \left(\frac{\nu_1 b_{2n+1,1}}{\sigma_{2n+1}} + r_{2n+1}\right) P_{2n+1}^{(\a2)}(x),
\end{eqnarray*}
where $r_\ell$ becomes negligible as $\ell$ increases. Therefore, if $n$ is sufficiently large then from (\ref{amm1}) we get
\begin{eqnarray} \label{casp0}
c_{2n} &\approx&   \frac{\nu_0 b_{2n,0}- \nu_2 b_{2n,2}}{\mu_{2n}\sigma_{2n}} =  \frac{\nu_0 b_{2n,0}- \nu_2 b_{2n,2}}{a_{2n}} \equiv \hat{c}_{2n},\\
c_{2n+1} &\approx& \frac{\nu_1 b_{2n+1,1}}{\mu_{2n+1} \sigma_{2n+1}} = \frac{\nu_1 b_{2n+1,1}}{a_{2n+1}} \equiv \hat{c}_{2n+1}. \label{casd0}
\end{eqnarray}
The statement follows from (\ref{bmn_as}) by observing that $a_\ell= O(\ell^{-1}).$ \qed\\

\begin{figure}[t]
\begin{center}
       \includegraphics[width=12.5cm,height=9cm]{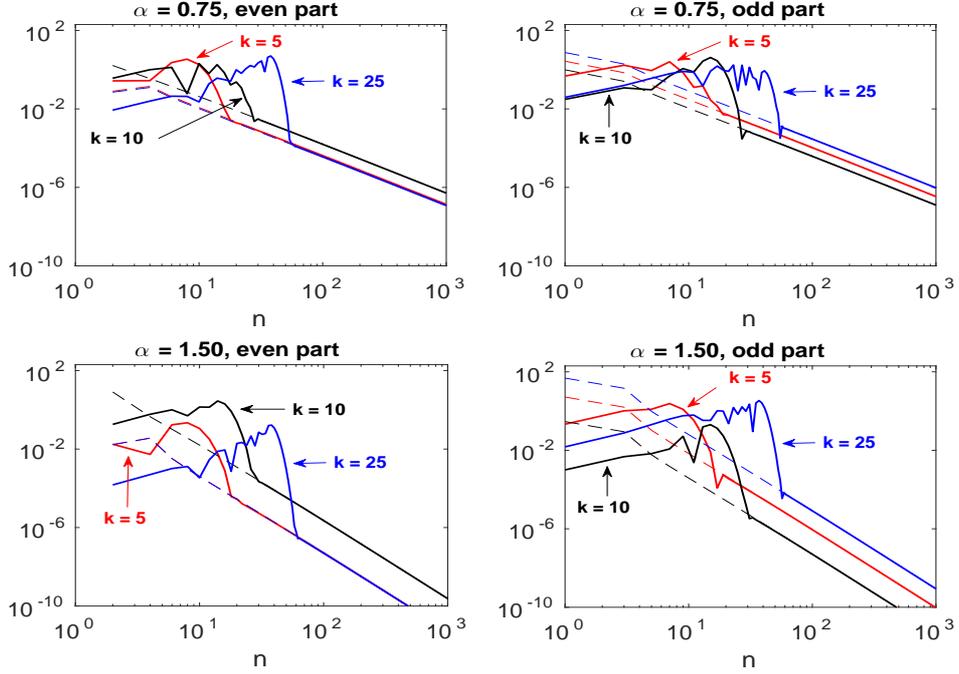}
\end{center}       
\caption{The coefficients  $|c_n|$ (solid lines) and $|\hat{c}_n|$ (dashed lines) for $q(x) = \mathrm{e}^x.$ }
\label{fig:asintc}
\end{figure}

In Figure~\ref{fig:asintc}, some numerical illustrations of the estimates in (\ref{casp0}) (left plots) and (\ref{casd0}) (right plots) 
have been reported for $q(x) = \mathrm{e}^x,$ $\alpha=0.75,1.5$ and $k=5,10,25.$ As one can see, for this example,  
such approximations of the coefficients are rather accurate for each $n\ge 3k.$\\

Before proceeding, we need the following notation. For each $\nu >0$ and $x\in (0,1)$ let 
\begin{equation}\label{J}
J(x,\nu,\mu) \equiv \int_0^x t^{\nu-1} (1-t)^{\mu-1} dt = \frac{x^\nu\,(1-x)^\mu}{\nu} {}_2F_1(1,\nu+\mu; 1+\nu;x) 
\end{equation}
where the second equality is an application of Euler's hypergeometric transformations \cite{Abram}.
By using a suitable Taylor expansion one deduces that
\begin{itemize}
 \item if $x \rightarrow 0^+$ then
 \begin{equation} \label{J0}
  J(x,\nu,\mu) = \frac{x^\nu\,(1-x)^\mu}{\nu} \left( 1 + O(x) \right);
 \end{equation}
\item if $\mu <0$ then ${}_2F_1(1,\nu+\mu; 1+\nu;1)=-\nu/\mu$ and consequently
\begin{equation} \label{J1}
J(x,\nu,\mu) = - \frac{x^\nu\,(1-x)^\mu}{\mu} \left( 1 + o(1) \right), \quad {\mbox as \,\,} x \nearrow 1^{-}.
\end{equation}
More precisely if $\mu<-1$ then $\left.\frac{d}{dx}{}_2F_1(1,\nu+\mu; 1+\nu;x)\right|_{x=1}=\frac{(\nu+\mu)\nu}{(\mu+1)\mu}.$
This implies that
\begin{equation} \label{J2}
J(x,\nu,\mu) = - \frac{x^\nu\,(1-x)^\mu}{\mu} \left( 1 + \frac{\nu+\mu}{\mu+1} (1-x) + o((1-x)) \right).
\end{equation}
\end{itemize}
We are now ready for the analysis of the entries of the local truncation error in (\ref{tau_m2}).

\begin{theorem}\label{Imn_teo}
Let $\al\in (0,2)$ and suppose that $q$ verifies the assumptions {\bf H1} and {\bf H2}.
If $k$ is the index of the eigenvalue and $N_0$ is sufficiently larger than $k,$ then for each $N\ge N_0$ one has
\begin{equation}\label{taumup}
|\tau_m^{(N)}| \leq \frac{C_0}{2\alpha N^{4\alpha+2}} \left(1 -\frac{m^2}{N^2}\right)^{-\alpha}, \qquad m=0,1,\ldots,N-1,
\end{equation}
where $C_0$ is a constant independent of $N.$
\end{theorem}
\underline{Proof} 
 By using an integral estimate, from (\ref{tau_m2}), (\ref{bmn_as}), (\ref{qmn_as}) and (\ref{cn_as}) one deduces 
that there exist a constant $C_0,$ independent of $N$ and $m,$ such that
\begin{equation}\label{Imn}
|\tau_m^{(N)}| \leq \sum_{n=N}^\infty |u_{mn}| |c_n| \leq C_0 \int_{N}^\infty (n^2-m^2)^{-\al-1}n^{-2\al-1} dn \equiv C_0 I_m^{(N)}.
\end{equation}
If $m>0$ and if  we apply the change of variable $n=mt^{-1/2}$ then, after some computation,  from (\ref{J})
with $\nu = 2\al+1$ and $\mu =-\al,$ we get
\begin{eqnarray*}
I_m^{(N)} &=& \frac{1}{2m^{4\al+2}} \,J \left(m^2/N^2, 2\al+1,-\al\right) \\
&=& \frac{1}{2(2\al+1) N^{4\al+2}} \,\left(1-\frac{m^2}{N^2}\right)^{-\al}\,{}_2F_1(1,\al+1; 2\al+2;m^2/N^2).
\end{eqnarray*}
It is not difficult to verify that this equality holds true also for $m=0.$
We observe that ${}_2F_1(1,\al+1; 2\al+2;x)$ is increasing over $[0,1]$ and its value at $x=1$ is $(2\al+1)/\al.$
Therefore
$$I_m^{(N)} \leq \frac{1}{2\al N^{4\al+2}} \,\left(1-\frac{m^2}{N^2}\right)^{-\al}$$
which, together with (\ref{Imn}), gives (\ref{taumup}). \qed\\

The immediate consequence of this result is that the first and the last entries of $\bftau^{(N)}$ behave like 
$O(N^{-4\al-2})$ and $O(N^{-3\al-2}),$ respectively. \\
The following theorem completes the error analysis in the eigenvalue approximation.
\begin{theorem}\label{teoconv}
If the hypotheses of the previous theorem hold true then there exist a constant $C,$ independent of $N,$
such that
\begin{equation}\label{ordine}
|\lambda-\lambda^{(N)}| \leq C N^{-(4\al+2)}.
\end{equation}
\end{theorem}
\underline{Proof} 
 From (\ref{errlam})-(\ref{limden}), we deduce that there exist a constant $C_1$ such that
$$|\lambda-\lambda^{(N)}| = |\bfc_N^T \bftau^{(N)} + (\bfxi_N-\bfc_N)^T \bftau^{(N)}|/|(\bfxi_N^T B_N \bfc_N| 
\leq C_1 |\bfc_N^T \bftau^{(N)}|.$$
Now, by virtue of Theorem~\ref{thcnas} we have that  $|c_m|\sim m^{-2\al-1}$ for each $m\ge S =S(k),$
being $k$ the index of the eigenvalue. 
Therefore, from the previous theorem, we get
\begin{eqnarray*}
|\bfc_N^T \bftau^{(N)}| &\leq&  \sum_{m=0}^{S-1} |c_m| |\tau_m^{N}| +
\sum_{m=S}^{N-1} |c_m| |\tau_m^{N}| \\
&\leq& \frac{C_2}{N^{4\al+2}} \left(1 + \int_{S}^{N-1} m^{-2\al-1} \left(1-\frac{m^2}{N^2}\right)^{-\al}dm \right),
\end{eqnarray*}
where $C_2$ is a further suitable constant. If we apply the change of variable $t=1-m^2/N^2$ we obtain
\begin{eqnarray}\nonumber 
\int_{S}^{N-1} m^{-2\al-1} \left(1-\frac{m^2}{N^2}\right)^{-\al}dm = 
\frac{1}{2 N^{2\al}} 
\int_{x_0}^{x_1} t^{-\al}(1-t)^{-\al-1} dt \equiv \frac{\tilde{I}^{(N)}}{2N^{2\al}},
\end{eqnarray}
where
\begin{equation}\label{x0x1}
x_0=(2N-1)/N^2\approx 0, \qquad x_1 = 1-S^2/N^2 \approx 1.
\end{equation}
It follows that (\ref{ordine}) holds true provided that $\tilde{I}^{(N)} = O(N^{2\al}).$
This fact can be verified, after some computation, for $\al=1.$ 
If $\al \in (0,1)$ then from (\ref{J}), with $\nu = 1-\al >0$ and $\mu=-\al,$ 
(\ref{J0})-(\ref{J1}) and (\ref{x0x1}) we get
\begin{eqnarray}
\nonumber \tilde{I}^{(N)} = 
                 J\left(x_1, 1-\al,-\al\right) -J\left(x_0, 1-\al,-\al\right) = O(N^{2\al}).
\end{eqnarray}
Finally, if $\al \in (1,2)$ then with an integration by parts and by using 
(\ref{J}) with $\nu = 2-\al >0$ and $\mu=-\al-1,$ see also (\ref{J0})-(\ref{J2}) and (\ref{x0x1}), we obtain
\begin{eqnarray*}
\tilde{I}^{(N)} &=& \frac{1}{1-\al} \left.t^{1-\al} (1-t)^{-\al-1} \right\vert_{x_0}^{x_1}
 -\frac{\al+1}{1-\al} \int_{x_0}^{x_1} t^{1-\al} (1-t)^{-\al-2} dt\\
 &=& \frac{1}{1-\al} \left.t^{1-\al} (1-t)^{-\al-1} \right\vert_{x_0}^{x_1} \\
 &-& \frac{1+\al}{1-\al} \left( J(x_1,2-\al,-1-\al) -J(x_0,2-\al,-1-\al)\right)\\
 &=& \frac{x_0^{1-\al} (1-x_0)^{-\al-1}}{\al-1} \left(1 + O(x_0)\right)\\
 &+& \frac{x_1^{1-\al} (1-x_1)^{-\al-1}}{1-\al} \left( 1-x_1(1-(1-2\al)(1-x_1)/\al + o(1-x_1))\right)\\
 &=& O (N^{\al-1}) + \frac{x_1^{1-\al} (1-x_1)^{-\al}}{1-\al} \left( 1+(1-2\al)x_1/\al + o(1))\right)\\
 &=& O (N^{\al-1}) + O(N^{2\al}) = O(N^{2\al}),
\end{eqnarray*}
which completes the proof. \qed\\ 

\section{Conditioning analysis}\label{condizio}
We now discuss the conditioning of the numerical eigenvalues with respect to a perturbation of the 
potential. For a fixed $N,$ let $\lambda^{(N)}\equiv \lambda_k^{(N)}(q)$ and $\bfxi_N \equiv \bfxi_{k,N}(q)$
be the $k$-th numerical eigenvalue and the corresponding 
eigenvector as defined in (\ref{genmat}). In addition, let
$y^{(N)}(x)=y^{(N)}_k(x)$ be the resulting approximation of $y(x)=y_k(x)$  specified in (\ref{yN}).
If we apply the matrix method to problem (\ref{fslp})-(\ref{bc}) 
with perturbed  potential $\hat{q}(x)\approx q(x)$ then (\ref{genmat}) becomes
$$ \left(A_N + \hat{Q}_N\right) \hat{\bfxi}_{k,N}(\hat{q}) = \lambda_k^{(N)}(\hat{q}) B_N \hat{\bfxi}_{k,N}(\hat{q})$$
where the entries of $\hat{Q}_N$ are given by, see (\ref{qmn}),
$$ \hat{q}_{mn} = \langle P_m^{(\a2)}, \hat{q} P_n^{(\a2)}\rangle_{\al}, \quad m,n=0,\ldots N-1.$$
With these notations, we can prove the following result which concerns the case of a small perturbation of $q,$
i.e. of the initial datum.

\begin{proposition}\label{Pcond}
If  $\|q-\hat{q}\|_\infty$ is small enough then 
$$\left|\lambda_k^{(N)} (q) -\lambda_k^{(N)} (\hat{q})\right| \lesssim \|q-\hat{q}\|_\infty, \quad k=0,1,\ldots,N-1.$$
\end{proposition}
\underline{Proof} 
We observe that $(A_N+Q_N,B_N)$ is a symmetric (hermitian) definite pair, see \cite{Stew} and Remark~\ref{Bsdp}.
Hence, by using standard arguments,                 
it is not difficult to obtain that if $\|q-\hat{q}\|_\infty$ is sufficiently small
then 
\begin{eqnarray*} \left|\lambda_k^{(N)} (q) -\lambda_k^{(N)}(\hat{q})\right| &=&
\left|\frac{\bfxi_{k,N}^T(q) \left(Q_N-\hat{Q}_N\right)\bfxi_{k,N}(\hat{q})}
{\bfxi_{k,N}^T(q) B_N\bfxi_{k,N}(\hat{q})}\right|\\
&\approx&
\left|\frac{\bfxi_{k,N}^T(q) \left(Q_N-\hat{Q}_N\right)\bfxi_{k,N}(q)}
{\bfxi_{k,N}^T(q) B_N\bfxi_{k,N}(q)} \right| \\
&=&\left|\frac{\displaystyle\int_{-1}^1 (q(x)-\hat{q}(x) ) \left(y_k^{(N)}(x)\right)^2\,dx}
 {\displaystyle\int_{-1}^1 \left(y_k^{(N)}(x)\right)^2\,dx} \right|\leq \|q-\hat{q}\|_\infty. 
\end{eqnarray*}
\qed\\

The numerical eigenvalues are therefore definitely well-conditioned with respect to a 
perturbation of the potential. This fact and Remark~\ref{qpoli} suggest to consider a perturbed problem with 
$q$ replaced by $q_L \in \Pi_L,$ the space of polynomials of maximum degree $L,$
where $L$ represents a further parameter to be specified by the user (in addition to $N,$ the
order of the generalized matrix eigenvalue problem). In this way, in fact, the computation of $Q_N$
is simple and the resulting approximation of the eigenvalue verifies

$$\left|\lambda_k -\lambda_k^{(N)}(q_L)\right| \leq 
\left|\lambda_k -\lambda_k^{(N)}(q)\right| + \left|\lambda_k^{(N)}(q) -\lambda_k^{(N)}(q_L)\right|$$
which is approximately equal to $\left|\lambda_k -\lambda_k^{(N)}(q)\right|$ if $q_L$ is chosen 
properly.\\
In more details, even though different strategies are possible, we decided to select $q_L$ as the following
partial sum of the Fourier-Legendre series of $q$
\begin{equation}\label{qL}
q_L (x) = \sum_{j=0}^L \tilde{q}_j P_j^{(0)} (x), \qquad \tilde{q}_j = \frac{\langle q, P_j^{(0)}\rangle_0}
{\langle P_j^{(0)}, P_j^{(0)}\rangle_0} = \frac{ (2j+1)\,\langle q, P_j^{(0)}\rangle_0}{2},
\end{equation}
which converges in uniform norm versus $q$ as $L$ increases if $q$ is analytic on the Bernstein ellipse.
One motivation that lead us to consider this approach is that $q$ and $q_L$ have the same  mean value 
(see (\ref{lamasq}) and the third example in Section~\ref{example}). It must be said, in fact, that 
all our experiments suggest that the upper bound  established in Proposition~\ref{Pcond} may be not 
sharp if $\bar{q}=\bar{\hat{q}},$ $N$ is sufficiently large and $k/N<1/2.$ 
Finally, we observe that we can approximate
$ \langle q, P_j^{(0)}\rangle_0 = \int_{-1}^1 q(x) P_j^{(0)}(x) dx$
by applying a standard Gaussian
quadrature formula or, for example, the algorithm described in \cite{Iser}.

\section{Numerical experiments}\label{example}
The described method has been implemented in \texttt{MATLAB} (version \texttt{R2015b}) and 
the generalized matrix eigenvalue problem (\ref{genmat}) has been solved by using the \texttt{eig} or the
\texttt{eigs} (with option \texttt{SM}) commands depending on the number of eigenvalues 
we were interested in. For the approximation of the coefficients $\tilde{q}_j$ in (\ref{qL}) 
we have applied the standard Gaussian quadrature formula with degree of precision $\max (2L+1,11)$  
(the function
\texttt{legpts} included in \texttt{Chebfun v4.3.2987} has been used for the computation of its nodes and weights).
Concerning the choice of $L,$ we always select it in such a way that $\|q-q_L\|_\infty$ is of the 
order of the machine precision.
For the estimate of the error in the approximation of $\lambda_k$ we consider as ``exact'' 
the corresponding numerical eigenvalue obtained with matrices of size $N_{\mathrm{true}} \gg N$ with $N_{\mathrm{true}} = 
N_{\mathrm{true}}(k,\alpha).$
In fact, it is important to underline the fact that if $\alpha \in (0,2)$ then the exact eigenvalues 
are not known in closed form  even for $q$ identically zero in $(-1,1)$ and that, differently with respect to
classical Sturm-Liouville problems, nowadays it is not yet available a 
well-established numerical software for the solution of the eigenvalue problem we have
studied in this paper.\\

\begin{figure}
\begin{center}
 \includegraphics[width=12.5cm, height = 9cm]{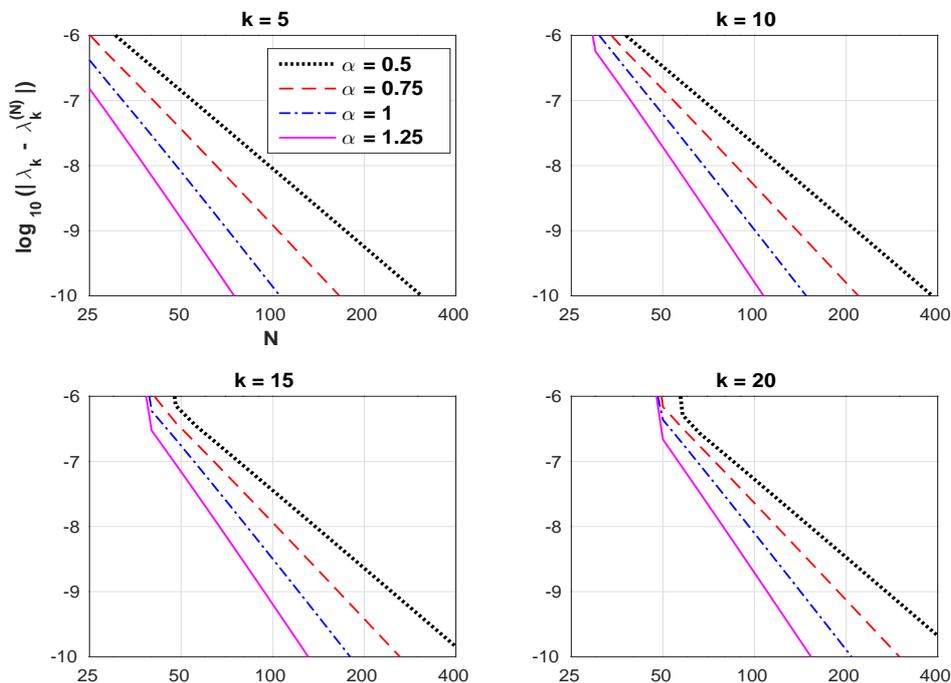}
\caption{Errors in the approximation of the eigenvalues for $q(x) = 2x(x+1).$}
 \label{figpoli}
\end{center}
\end{figure}

\begin{figure}
\begin{center}
 \includegraphics[width=12.5cm, height = 9cm]{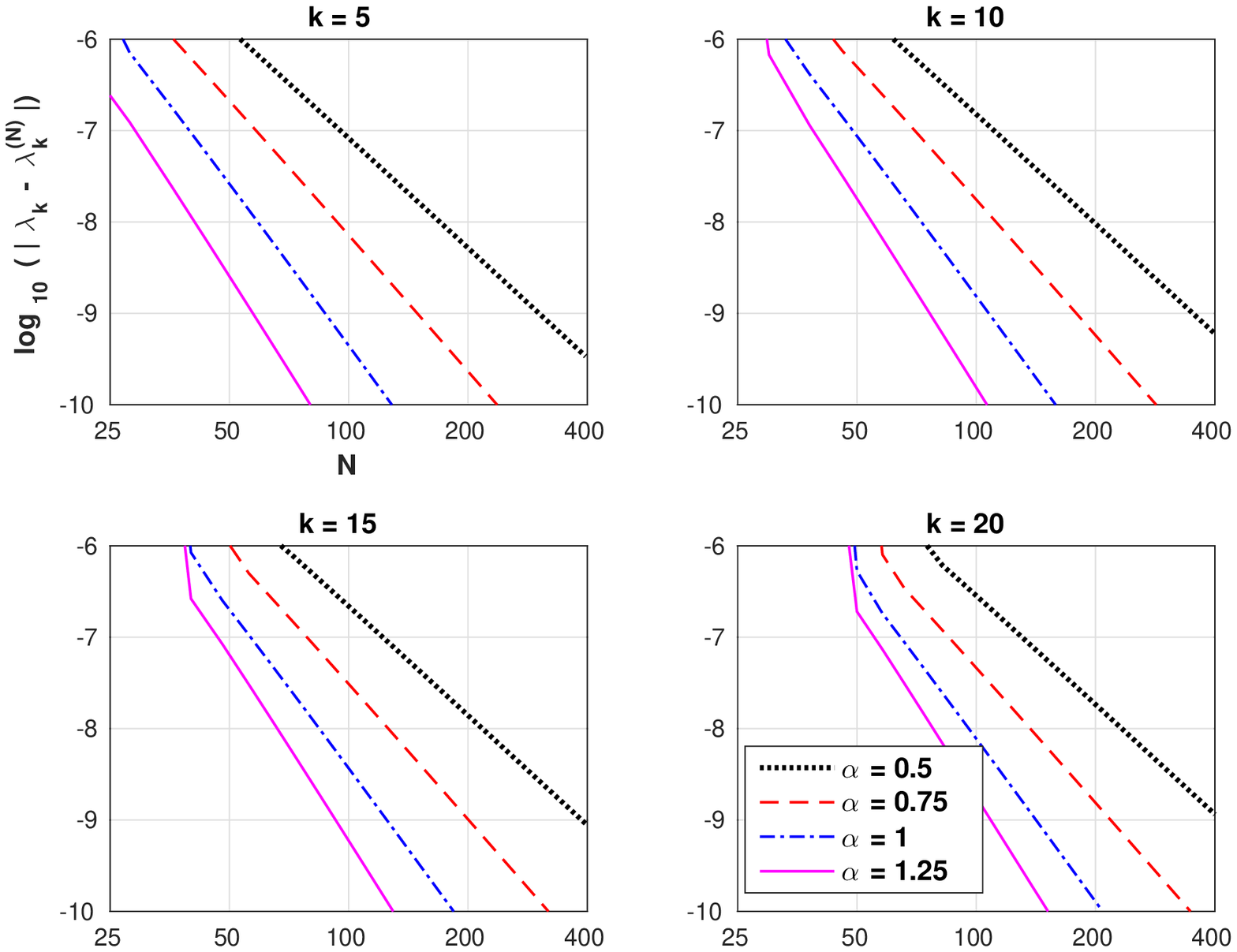}
\caption{Errors in the approximation of the eigenvalues for $q(x) = \mbox{e}^{\pi(x+1)/2}.$}
 \label{figexp}
\end{center}
\end{figure}

As first examples, we solved the problems with $q(x) =2x(x+1)$ and $q(x) = \mbox{e}^{\pi(x+1)/2}$ 
for several values of $\alpha$ and $N.$ Clearly, see (\ref{qL}),  for the first potential we set $L=2$ while 
for the second one we choose $L=15.$  In Figures~\ref{figpoli},\ref{figexp}, the errors in the resulting
approximation of the eigenvalues have been reported with, as we are going to do for all examples,
a logarithmic scale on the abscissae.
As one can see, in both cases, as soon as $N$ is sufficiently larger than 
the index of the eigenvalue $k,$ it results
$$ \log_{10}(|\lambda_k -\lambda_k^{(N)}|) \sim  -p \log_{10}(N).$$
Estimates of the various values of $p,$ determined with a least-square approximation, are listed in
Table~\ref{taborder}. It is evident that $p \approx (2+4\alpha)$ in agreement with Theorem~\ref{teoconv}.\\

\begin{table}
\begin{center}
 \begin{tabular}{|c||c|c|c|c|}
 \hline
  \multicolumn{5}{|c|}{~} \vspace{-.25cm}\\
   \multicolumn{5}{|c|}{$q(x) = 2x(x+1)$} \\
   \hline
   $k$ & $\alpha = 0.50$ & $\alpha = 0.75$ & $\alpha = 1.00$ & $\alpha = 1.25$ \\
   \hline
5 & 3.99 & 4.94 & 5.86 & 6.79 \\
10 & 3.99 & 4.96 & 5.94 & 6.92 \\
15 & 3.99 & 4.94 & 5.89 & 6.84 \\
20 & 4.00 & 4.96 & 5.94 & 6.92 \\
\hline
  \multicolumn{5}{|c|}{~} \vspace{-.25cm}\\
   \multicolumn{5}{|c|}{$q(x) = \mbox{e}^{\pi(x+1)/2}$} \\
   \hline
   $k$ & $\alpha = 0.50$ & $\alpha = 0.75$ & $\alpha = 1.00$ & $\alpha = 1.25$ \\
   \hline
5 & 4.00 & 4.96 & 5.91 & 6.85 \\
10 & 4.00 & 4.95 & 5.91 & 6.88 \\
15 & 4.01 & 4.95 & 5.90 & 6.85 \\
20 & 4.01 & 4.95 & 5.92 & 6.91 \\
\hline
 \end{tabular}
\caption{Estimates of the order of convergence in the eigenvalue approximations.}\label{taborder}
\end{center}
\end{table}

\begin{figure}[h]
\begin{minipage}[!t]{0.9\textwidth}
 \includegraphics[width=.9\textwidth, height = 6.5cm]{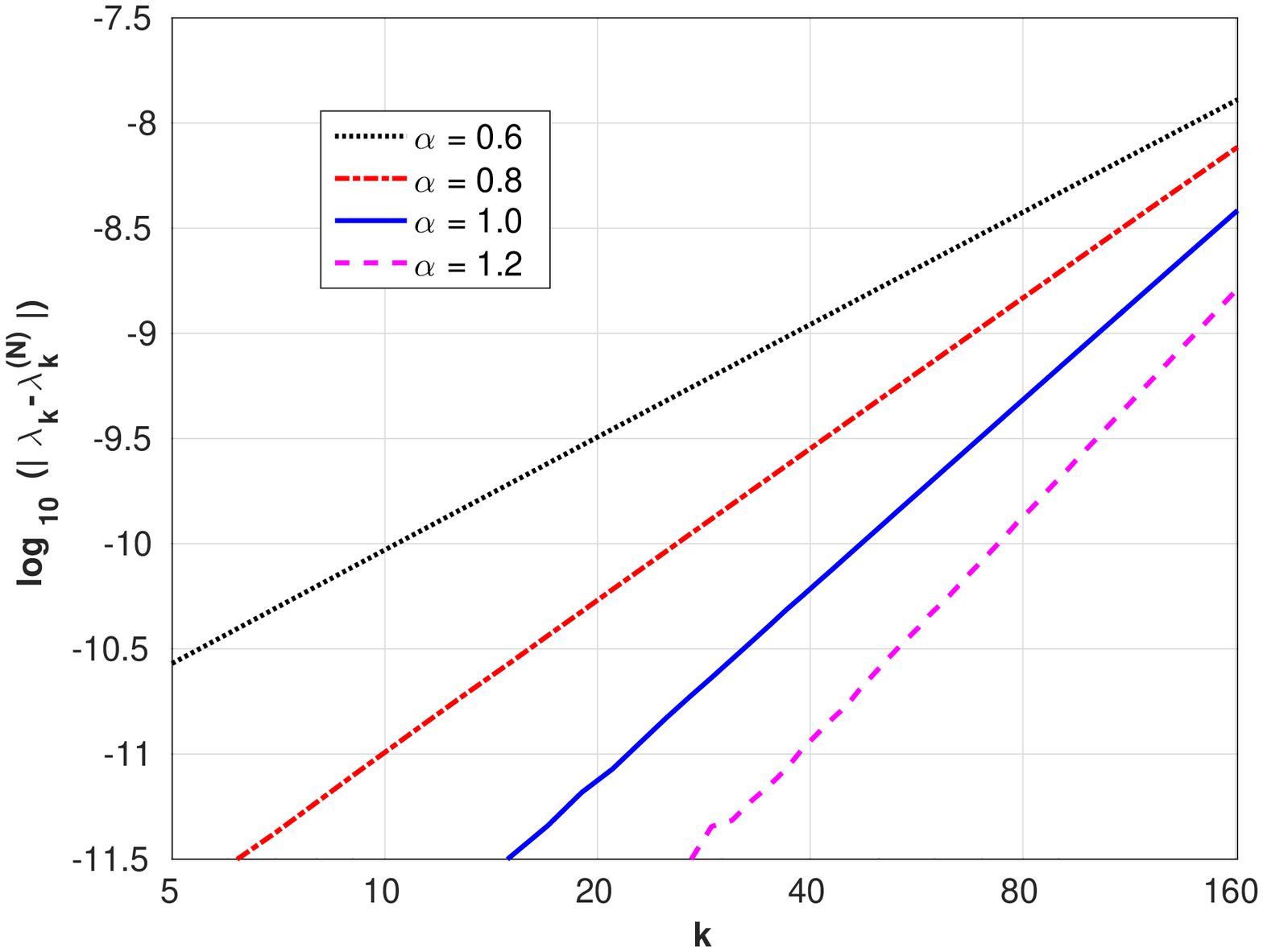}
\end{minipage}
\hspace{-0.1\textwidth}
\begin{minipage}[!t]{0.1\textwidth}
\begin{tabular}{|c|c|}
\hline
$\alpha$ & $r$\\
\hline
0.6  & 1.82 \\
0.8  & 2.41 \\
1.0  & 3.01 \\
1.2  & 3.60 \\
\hline
 \end{tabular}
\end{minipage}
\caption{Errors in the approximation of the eigenvalues versus their index for 
$q(x) = (x+1)/(2(x^2+1))$ with $N=320$ and $L=37.$}
 \label{figtrig}
\end{figure}

The aim of this second example is that of showing experimentally the growth of the error
in the eigenvalue approximation with respect to the index $k$ for a fixed $N.$ In all our 
experiments it seems that if $q$ verifies the assumptions {\bf H1} and {\bf H2} and 
if $k/N$ is sufficiently small then 
$$|\lambda_k -\lambda_k^{(N)}| = O\left(k^{r} N^{-2-4\al}\right), \qquad r = 3\al.$$
This can be partially explained by observing that $\lambda_k$ grows like $k^\al,$ 
see (\ref{lamasq0})-(\ref{lamasq}), and, 
consequently, see (\ref{tau_m2}), (\ref{nu0a})-(\ref{nu1b}) and (\ref{casp0})-(\ref{casd0}),
the first entries of the local truncation error behave like $k^{2\al}.$ Finally,
by using (\ref{normalcn})-(\ref{limden}) and (\ref{yasint}), one deduces that
the denominator in (\ref{errlam}) decreases like $k^{-\al}$ (this fact is rather
simple to be verified for $\al=2$ and $q\equiv 0$ over $(-1,1)$ since the exact eigenfunctions are
known in closed form). All these arguments lead to the value of $r$ previously specified.
As an example, in Figure~\ref{figtrig} we report the errors in the eigenvalue
approximations versus the index for $q(x) = (x+1)/(2(x^2+1))$
and for four values of $\al.$ In addition, we list estimates of the corresponding $r$'s 
determined with a least-square fitting.\\

\begin{figure}
\begin{minipage}[!t]{0.9\textwidth}
 \includegraphics[width=.9\textwidth, height = 6.5cm]{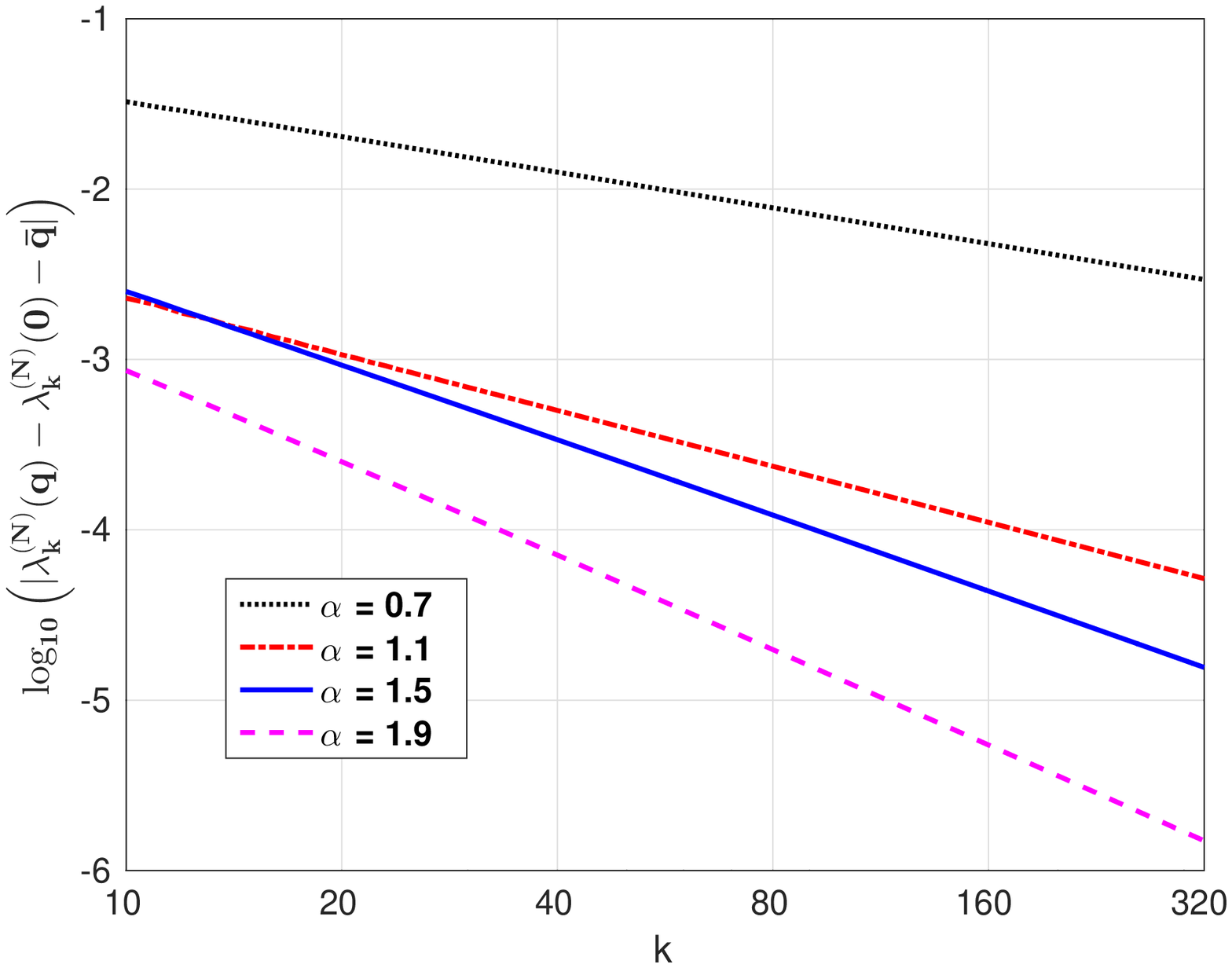}
\end{minipage}
\hspace{-0.1\textwidth}
\begin{minipage}[!t]{0.1\textwidth}
\begin{tabular}{|c|c|}
\hline
$\alpha$ & $\eta$\\
\hline
0.7  & 0.70 \\
1.1  & 1.09 \\
1.5  & 1.48 \\
1.9  & 1.86 \\
\hline
 \end{tabular}
\end{minipage}
\caption{Errors of the asymptotic estimate $\lambda_k^{(N)}(q) \approx \lambda_k^{(N)}(0) + \bar{q}$ for 
$q(x) = -\cos(3 x) + \sin(2 x)$ with $L = 18$ and $N = 1000.$}
 \label{figmedia}
\end{figure}

The third example is in support of the asymptotic estimate in (\ref{lamasq}). In particular, 
in Figure~\ref{figmedia}, we report $\log_{10}\left(|\lambda_k^{(N)}(q) -\lambda_k^{(N)}(0) -\bar{q}|\right)$
versus $k,$ for $q(x) =  -\cos(3 x) + \sin(2 x)$, with $L = 18$, $N = 1000$, and
$\alpha = 0.7,1.1,1.5,1.9.$ As one can see, for each $\alpha,$  
$\lambda_k^{(N)}$ approaches $\lambda_k^{(N)}(0) +\bar{q}$ as $k$ increases and this is in perfect 
agreement with (\ref{lamasq}) by considering also that $q$ and $q_L$ in (\ref{qL}) have the same mean value. 
As done in the previous examples, we apply a least-square fitting
to determine the values of $\eta$ such that $|\lambda_k^{(N)}(q) -\lambda_k^{(N)}(0) -\bar{q}| = O(k^{-\eta})$
and the resulting exponents are listed in the table on the right of the same figure. For this example,
we observe that $\eta \approx \alpha.$ \\

\begin{figure}[t]
\begin{center}
 \includegraphics[width=12.5cm, height = 9cm]{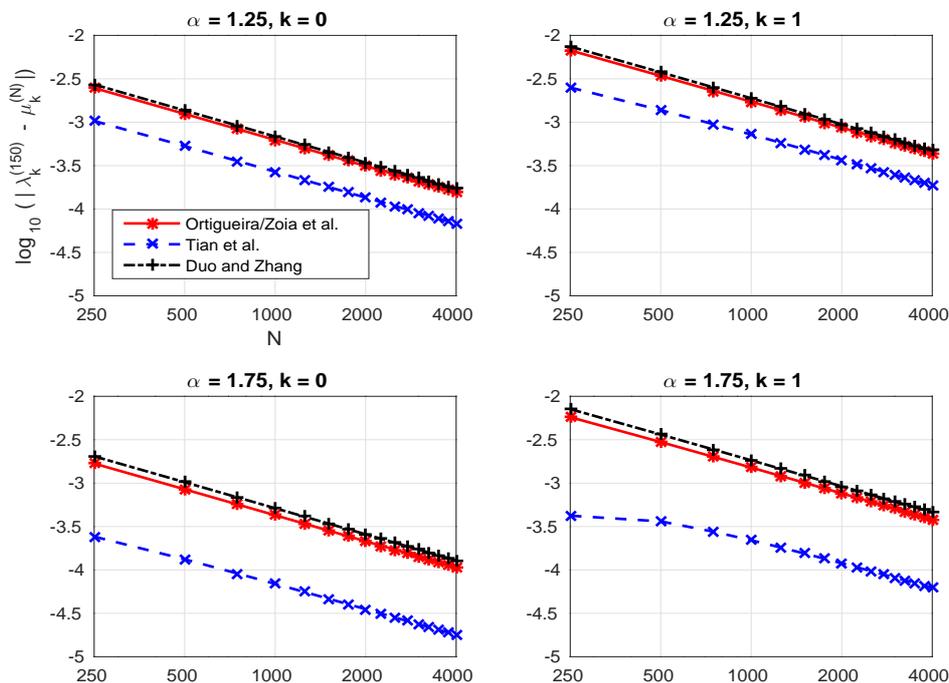}
\caption{Comparison of the estimates of the first two eigenvalues of the 
 infinite potential well problem provided by our method with $N=150$ and by the schemes
 proposed in \cite{Or06,ZRK07,TZD15,DZ15} with various values of $N.$}
 \label{confronti}
\end{center}
\end{figure}

Finally, we consider the infinite potential well, i.e. $q(x)=0$ for each $x\in (-1,1),$ with 
$\alpha=1.25$ and $\alpha = 1.75.$
The variations in the numerical approximations of its
first two eigenvalues (sometimes called the energies of the ground and of the first excited states) provided by the method 
proposed in this paper are of the order of the machine precision for each $N \ge 150.$
We then compare $\lambda_0^{(150)}$ and $\lambda_1^{(150)}$ 
with the numerical eigenvalues given by the methods proposed by 
\begin{itemize}
 \item Ortigueira/Zoia et al. in \cite{Or06,ZRK07};
 \item Tian et al. in \cite{TZD15};
 \item Duo and Zhang in \cite{DZ15}.
\end{itemize}
We shall call $\mu_k^{(N)}$ the estimate of
$\lambda_k$ provided by one of the previous three methods with a matrix of order $N.$
The results so obtained are reported in Figure~\ref{confronti}.
It is evident that $|\lambda_k^{(150)}-\mu_k^{(N)}|$ always decreases at the same rate; more
precisely we have verified that such difference behaves like $O(N^{-1}).$
Finally, it is important to mention the fact  that we have done  similar 
experiments with other potentials and that  the method that we propose turns out to 
be absolutely competitive with the other three ones in all our tests.

\end{document}